\documentclass{article}
\usepackage{color}
\usepackage{amsmath}
\usepackage{amscd}
\usepackage{latexsym}
\usepackage{amssymb}
\usepackage[dvips]{graphicx}
\usepackage{amsfonts}
\usepackage{amssymb}
\usepackage[utopia]{mathdesign}
\usepackage{mathtools}
\usepackage[version=3]{mhchem}
\usepackage{gensymb}
\usepackage{color}
\def\trd{\textcolor{red}}
\def\tpr{\textcolor{purple}}
\def\tgr{\textcolor{green}}
\def\tor{\textcolor{orange}}
\title{Cohomology of Lie semidirect products\\ and poset algebras}
\author{Vincent E. Coll, Jr. and Murray Gerstenhaber}

\begin{document}
\maketitle

\newtheorem{theorem}{Theorem}
\newtheorem{corollary}{Corollary}
\newtheorem{lemma}{Lemma}
\def\trd{\textcolor{red}}
\def\tpr{\textcolor{purple}}
\def\tgr{\textcolor{green}}
\def\tor{\textcolor{orange}}
\renewcommand{\abstractname}\!\!{}
\newcommand{\C}{\ensuremath{\mathbb{C}}}
\newcommand{\R}{\ensuremath{\mathbb{R}}}
\newcommand{\Z}{\ensuremath{\mathbb{Z}}}
\newcommand{\calF}{\ensuremath{\mathcal{F}}}
\newcommand{\calP}{\ensuremath{\mathcal{P}}}
\newcommand{\pr}{\ensuremath{\preceq}}
\newcommand{\op}{\ensuremath{\mathrm{op}}}
\newcommand{\ctr}{\mathfrak{c}}
\newcommand{\g}{\mathfrak{g}}
\newcommand{\h}{\mathfrak{h}}
\newcommand{\frk}{\mathfrak{k}}
\newcommand{\G}{\mathfrak{G}}
\newcommand{\n}{\mathfrak{n}}
\newcommand{\bfk}{\mathbf{k}}
\newcommand{\w}{\mathbf{w}}
\newcommand{\bfl}{\mathbf{\lambda}}
\newcommand{\m}{\mathfrak{m}}
\newcommand{\mbfz}{\mathbf{0}}
\renewcommand{\b}{\mathfrak{b}}
\newcommand{\ad}{\ensuremath\operatorname{ad}}
\newcommand{\diag}{\ensuremath\operatorname{diag}}
\newcommand{\Hom}{\ensuremath\operatorname{Hom}}
\newcommand{\Sk}{\ensuremath\operatorname{Sk}}
\newcommand{\myrightleftarrows}[1]{\mathrel{\substack{\xrightarrow{#1} \\[-.9ex] \xleftarrow{#1}}}}
\vspace{-2mm}
\date\!\!{}
{\noindent \textit{Department of Mathematics, Lehigh University, Bethlehem, PA 18015; Department of Mathematics, University of
Pennsylvania, Philadelphia, PA 19104-6395}
\vspace{-7mm}

\begin{abstract}\noindent  
When $\h$ is a toral subalgebra of a Lie algebra $\mathfrak g$ over a field $\mathbf k$, and $M$ a $\mathfrak g$-module on which $\mathfrak h$ also acts torally, the Hochschild-Serre filtration of the Chevalley-Eilenberg cochain complex admits a stronger form than for an arbitrary subalgebra. For a semidirect product $\mathfrak g = \mathfrak h \ltimes \mathfrak k$ with $\mathfrak h$ toral one has $H^*(\g, M) \cong \bigwedge\h^{\vee} \bigotimes H^*(\mathfrak k,M)^{\mathfrak h} = H^*(\mathfrak h, \mathbf k)\bigotimes H^*(\mathfrak k,M)^{\mathfrak h}$, and for a Lie poset algebra $\g$, that $H^*(\mathfrak g, \mathfrak g)$, which controls the deformations of $\mathfrak g$, can be computed from the nerve of the underlying poset. The deformation theory of Lie poset algebras, analogous to that of complex analytic manifolds for which it is a small model, is illustrated by examples.
\end{abstract}

\section{Introduction}
A Lie algebra $T$ acts torally on a vector space $V$ over $\bfk$ if all its elements act semisimply, or equivalently, can be brought into diagonal form over the algebraic closure $\bar \bfk$ of $\bfk$. The actions of its elements must  then commute,  cf.  \cite[p.34]{Humphreys:Lie}, so a toral subalgebra of a Lie algebra $\g$, i.e.,  one whose adjoint action is toral, is necessarily Abelian; it will here generally be denoted by $\h$. If  $\h$ also acts torally on a $\g$-module $M$, then it does so on the Chevalley-Eilenberg cochain complex $C^*(\g, M)$, and therefore on the cohomology $H^*(\g, M)$. Viviani's Lemma, \S2, implies that the Hochschild-Serre filtration of $C^*(\g,M)$ admits a stronger form than when $\h$ is an arbitrary subalgebra. When $\g$ is a semidirect product $\h\ltimes \frk$, this yields the following:
\begin{equation}\label{general}
H^*(\g, M) \; \cong\; \bigwedge\h^{\vee} \bigotimes H^*(\frk,M)^{\h} \,=\, H^*(\h,\frk)\bigotimes H^*(\frk,M)^{\h}.
\end{equation}
Here $\h^{\vee}$ is the dual space of $\h$ and 
$H^*(\frk,M)^{\h}$ consists of the invariants of $H^*(\frk,M)$ under the operation of $\h$; the last equation follows since $\h$ is Abelian.

Let $\cal P$ = $\{i,j, \dots \}$ be a finite poset with partial order $\preceq$. The associative poset algebra $A = A(\cal P, \bfk)$ is the span over $\bfk$ of elements $e_{ij}, \, i\preceq j$ with multiplication given by setting $e_{ij}e_{j'k}= e_{ik}$ if $j=j'$ and 0 otherwise. The trace of an element $\sum c_{ij}e_{ij}$ is $\sum c_{ii}$; the Lie poset algebra $\g = \g(\cal P, \bfk)$ is the Lie subalgebra  of $A$ of all elements of trace zero. If $\#\cal{P}$ = $N$, in which case we may assume that the underlying set of $\cal P$ is the set of integers $\{1,\dots, N\}$ with partial order compatible with the linear order, then $A$ and $\g$ may be viewed as subalgebras of the algebra of all upper triangular $N \times N$ matrices over $\bfk$. Let $\b$ be the Borel subalgebra of this algebra consisting of matrices of trace zero and $\h$ its Cartan subalgebra of diagonal matrices.  Any subalgebra $\g$ with $\b \supset \g \supset \h$ is then a Lie poset algebra, a condition which may be used to define Lie poset algebras more generally. For $\g$ is then the span over the field $\bfk$ of $\h$ and of those $e_{ij}$ which it contains, and there is a partial order on $\{1,\dots,N\}$ compatible with the linear order defined by setting $i\prec j$ whenever $e_{ij} \in \g$.

Suppose now that the characteristic $p$ of $\bfk$ is either zero or greater than $N$. We compute the cohomology of $\g$ in two special cases, those where the module is $\bfk$ with trivial operation, and where it is $\g$ itself with adjoint operation, the latter being the case essential for deformation theory. Viewing $\cal P$ as a category, let $\Sigma$ be its nerve and $\Sigma^+$ be the complex obtained by adjoining a unique simplex of dimension -1 serving as boundary of every 0-simplex. Only the zero-dimensional cohomology is affected by this adjunction; one has $H^n(\Sigma^+, \bfk) = H^n(\Sigma,\bfk)$ for $n > 0$. We show that 
\begin{eqnarray}\label{intro_main}
 H^*(\frk,M)_{\mathbf 0}\;  \cong \;  H^*(\Sigma^+, \bfk),
 \end{eqnarray}
 where on the left one has Chevalley-Eilenberg cohomology and on the right simplicial cohomology with coefficients in $\bfk$. Combining this with (\ref{general}), one has
\begin{eqnarray*}
H^*(\g,\g)\; = \; \bigwedge \h^{\vee} \bigotimes H^*(\Sigma^+, \bfk) \; = \; H^*(\h,\bfk) \bigotimes H^*(\Sigma^+, \bfk).
\end{eqnarray*}
It follows that the space of infinitesimal deformations $H^2(\g,\g)$ of $\g$ is a direct sum of three components,
\begin{equation}\label{infinitesimals}
H^2(\g,\g) \;=\; (\bigwedge{}\!^2\,\h^{\vee} \bigotimes \ctr)\; \;\bigoplus\; \;(\h^{\vee} \bigotimes H^1(\Sigma,\bfk))\; \;\bigoplus\; \;H^2(\Sigma, \bfk).
\end{equation}
The first consists of infinitesimal deformations, necessarily in a non-commutative direction, of the structure of $\h$ alone. The second consists of infinitesimal deformations of the action of $\h$ on the $e_{ij}$.  The third consists of infinitesimal deformations of $\g$ determined by the complex $\Sigma$ alone, and in effect, only on the underlying topology of the geometric realization of its nerve.  In the final section we give simple examples of global deformations whose infinitesimals are of each of these types.

The decomposition \eqref{infinitesimals} is analogous to that of the space of infinitesimal deformations of a complex analytic manifold $\cal X$.   Those infinitesimal deformations of $\cal X$ to other complex manifolds were identified by by Fr\"olicher and Nijenhuis, \cite{FN:Stability}, with $H^1(\cal X, \cal T)$, where $\cal T$ denotes the sheaf of germs of holomorphic tangent vectors on $\cal X$.\footnote{This breakthrough  paved the way for the  work of Kodaira and Spencer, \cite {KS:Deformations}. Earlier,  Teichm\"uller,  \cite{Teich:quasikonforme}  had defined infinitesimal deformations of a Riemann surfaces, identifying these with its quadratic differentials, but his methods could not be extended to higher complex dimensions.} However, the full space of infinitesimal deformations of $\cal X$ has components not recognized until the introduction of algebraic  deformation theory, cf. \cite{GS:Monster}. When $\cal X$ is projective (and probably more generally), there is a single associative algebra $A$ built from $\cal X$ and a basic isomorphism, \cite{GS:Monster}, of cohomology rings,
$$H^*(A, A) \quad\cong\quad  H^*(\cal X, \bigwedge \cal T).$$
In particular,
\begin{equation}\label{complex}
H^2(A,A) \,\cong \,H^2({\cal X}, {\bigwedge} {\cal T}) \,
  =\, H^0 ({\cal X}, {\bigwedge}{}\!^2 {\cal T}) \bigoplus H^1({\cal X}, {\bigwedge}{}\!^1 {\cal T}) \bigoplus  H^2 ({\cal X}, {\bigwedge}{}\!^0 {\cal T}).
\end{equation}
The first component on the right consists of infinitesimal deformations of $\cal X$ to spaces whose function sheaves are sheaves of non-commutative rings,  cf. \cite[p. 250]{GS:Monster}. Such deformations are an essential aspect of  quantization; for an exposition and summary of the rich history of this idea, cf. \cite{DitoSternheimer:Genesis}. The second consists of the classical Fr\"olicher-Nijenhuis infinitesimals. The third component, which is just $H^2(\cal X, \C)$, has been called the ``mysterious" one by Kontsevich as it is difficult to see how elements of the second cohomology group of the underlying topological space of $\cal X$ can produce deformations. However, in the context of our present results, which provide a small model of the complex analytic case, we show how this happens.

By contrast with \eqref{intro_main}, in the associative case (where $\bfk$ may even be an arbitrary commutative unital coefficient ring), we have simply
$H^*(A,A) \; \cong\;  H^*(\Sigma, \bfk).$
This follows readily from the basic proposition that Hochschild cohomology can be computed relative to any separable subalgebra containing the unit element, cf. \cite{GS:SC=HC}, but no analog is known for Lie algebras.  Although the Euler-Poincar\'e characteristic of  an associative poset algebra can be arbitrary,  that of $H^*(\g,M)$ in \eqref{general} is zero  (like that of any finite-dimensional Lie algebra, cf. \cite{Goldberg:Lie}\footnote{Shown there only for the trivial module $\bfk$ but the proof holds for all finite-dimensional ones.For another proof which, by deformation theory extends to certain infinite dimensional algebras, cf. \cite{GG:WeylCohom}.}) because of the factor $\bigwedge\h^{\vee}$.
A third significant difference between the Lie and associative cases is that  associative poset algebras $A$ (over an arbitrary commutative unital ring $\bfk$) are quasi self-dual in the sense of \cite{G:Self-dual}, i.e. there is an isomorphism $H^*(A,A) \cong H^*(A, A^{\vee\,\op})$, where ``op'' denotes the interchange of left and right operations. It follows that  for such algebras $H^*(A,A)$ is a contravariant functor of $A$, but for Lie poset algebras $H^*(\g, \g)$ is generally not isomorphic to $H^*(\g, \g^{\vee})$. (Note that $\g =  \g^{\op}$ since every Lie algebra is isomorphic to its opposite: send every element to its negative.)

\section{Viviani's lemma}\label{Viviani}
Suppose that a Lie algebra $T$ acts torally on a vector space $V$ over $\bfk$; the action of $\tau \in T$ on $v \in V$ will then be denoted by $[\tau,v]$. When $\bfk$ is algebraically closed, the space $V$ splits into a direct sum of weight spaces, the weights $\w$ being elements of the dual space $T^{\vee}$; an element $v$ is in the weight space $V_{\w}$ if $\tau v = \w(\tau)v$ for all $\tau \in T$. However, even when $k$ is not algebraically closed, the weight $\mbfz$ space $V_{\mbfz}$, consisting of the  invariants under $T$ (elements it annihilates), is well-defined and a direct  summand of $V$, although it may be reduced to the zero element of $V$. 
If $T$ acts torally and compatibly on a Lie algebra $\g$ and $\g$-module $M$, i.e., if $[\tau,[g,g']] = [[\tau, g],g'] +[\g,[\tau, g']]$ and $[\tau, [g,m]] = [[\tau, g], m] + [g, [\tau, m]]$ for all $\tau \in T; g,g' \in \g$ and $m \in M$, then it also acts on the Chevalley-Eilenberg complex by setting 
$$[\tau, F](g_1,\dots, g_n)\;  = \; [\tau, F(g_1,\dots, g_n)]\;  - \; \sum_{i=1}^nF(g_1, \dots, [\tau, g_i],\dots,g_n),$$
where $g_1, \dots, g_n  \in \g$, and $F \in C^n(\g, M)$ is an $n$-cochain of $\g$ with coefficients in $M$. This action descends to the cohomology $H^*(\g, M)$ which, when $\bfk = \bar\bfk$, therefore also decomposes into a direct sum of weight spaces. Denoting the group of $n$-cochains of $\g$ with coefficients in $M$ by $C^n$, the Chevalley-Eilenberg coboundary operator $\delta :C^n \to C^{n+1}$ is defined\footnote{This definition, allowing the use of an arbitrary $\g$-module $M$, is probably due to Hochschild and seems to have appeared first in \cite{Hoch-Serre:Lie}. Chevalley and Eilenberg, \cite{Ch-E}, initially considered only cohomology with trivial coefficients, in which case the first term on the right does not appear} by setting
\begin{multline*}
(\delta F)(g_0, g_1, \dots, g_n)\quad=\quad
\sum_{0\le i\le n}(-1)^i[ g_i, F(g_0\dots,\hat g_i,\dots g_n)] \; +\\
\sum_{0\le i < j \le n}(-1)^{i+j}F([g_i,g_j],g_0\dots, \hat g_i,\dots, \hat g_j, \dots g_n),
\end{multline*}
where $\hat g$ indicates omission of the argument $g$.  Following \cite{Viviani:Melikian}, it is useful to rearrange the terms on the right by taking first all those in which $\ad g_0$ appears either as an operator on $M$ or on $\g$.  Denoting by  $\iota_{g_0} F$ the $(n-1)$-cochain defined by setting  $(\iota_{g_0})F(g_1,\dots,g_{n-1} )= F(g_0, g_1,\dots,g_{n-1})$, the coboundary then takes the form
\begin{multline*}
(\delta F)(g_0, g_1, \dots, g_n)\; =\\
[g_0, F(g_1,\dots,g_n)] \; + \sum_{1 \le i \le n}(-1)^i F([g_0,g_i],g_1,\dots, \hat g_i, \dots, g_n)\\
-\; (\delta (\iota_{g_0}\!\ F))(g_1,\dots, g_n).
\end{multline*}
Since $(\delta F)(g_0, g_1, \dots, g_n) =(\iota_{g_0}(\delta F))(g_1,\dots, g_n)$, this can be rewritten as
\begin{multline}\label{coboundary}
(\iota_{g_0}(\delta F))(g_1,\dots, g_n) \; +\; (\delta (\iota_{g_0}F))(g_1,\dots, g_n)=\\
[g_0, F(g_1,\dots,g_n)] \; + \sum_{1 \le i \le n}(-1)^i F([g_0,g_i],g_1,\dots, \hat g_i, \dots, g_n).
\end{multline}

\begin{theorem}\label{cochain map}
Suppose that a toral subalgebra $\h$ of $\g$ acts torally also on a $\g$-module $M$. If $F\in C^n(\g, M)$ is the homogeneous component of weight $\mathbf{w} \in \h^{\vee}$, then for all $h \in \h$ one has
\begin{equation}\label{cochain}
(\iota_h(\delta F))(g_1,\dots, g_n) \; +\; (\delta (\iota_h F))(g_1,\dots, g_n)\; = \; 
\mathbf{w}(h)F(g_1,\dots,g_n).
\end{equation}
In particular, if $F$ is a homogeneous cochain of weight zero then
\begin{equation}\label{cobound}
(\iota_h(\delta F))(g_1,\dots, g_n) \; \; =\; \;  -(\delta (\iota_h F))(g_1,\dots, g_n)
\end{equation}
for all $h \in \h$, i.e., $\iota_h$ is (up to sign) a cochain mapping when restricted to the subcomplex $C^*(\g,M)_{\mathbf 0}$ of cochains of weight $\mathbf{0}$.
\end{theorem}

\noindent \textsc{Proof.} Assume for the moment that $\bfk = \bar\bfk$. When the arguments are all homogeneous elements of $\g$ it is easy to see that \eqref{cochain} holds . This implies, however, that it is true in general, proving the first assertion. The second follows. $\Box$\medskip

\noindent  When $\bfk$ is extended to $\bar\bfk$ every cochain decomposes into the sum of its homogeneous parts under the operation of $\h$, but note that the part of weight $\mbfz$ is already defined over the original field $\bfk $; it is just the invariants under $\h$. 
\medskip

\noindent\textbf{Corollary (Viviani's Lemma)}
Any cocycle $F^n \in C^n(\g,M)$  is cohomologous to its homogeneous part of weight zero, whence
 \begin{equation}\label{weight zero}
 H^n(\g,M)\; \;  = \; \; H^n(\g,M)_{\mbfz} \;  =  \; H^n(\g,M)^{\h}. \end{equation}
\noindent \textsc{Proof.} Again, assume for the moment that $\bfk= \bar\bfk$.
If $F$ is a cocycle then so are all its homogeneous summands, so (\ref{cochain}) implies, for all $\w \in \h^{\vee}$, that if $F_{\w}$ is the homogeneous part of $F$ of weight $\w$, then $\w(h)F_{w}$ is a coboundary. If $\w \ne {\mbfz}$ then there is some $h \in \h$ such that $\w(h)\ne 0$, so $F_{\w}$ is a coboundary. Therefore, all that remains of the left side of (\ref{weight zero}) is the right side, whose definition does not require that $\bfk = \bar\bfk$.
$\Box$\smallskip

\noindent As  $H^n(\g,M)_{\mbfz} = H^n(\g,M)^{\h}$, we may use the notations interchangeably.\smallskip

Viviani's Lemma, \cite{Viviani:Melikian} (independently discovered but later, \cite{G:xLie}), does not exhaust for us the content of Theorem \ref{cochain map} as (\ref{cobound}) will also be essential in what follows.  The lemma applies in particular to semidirect products $\g = \h \ltimes \frk$ where $\h$ is toral
\footnote{The cohomology of Lie semidirect products $\g \ltimes \frk$, with $\g$ arbitrary and coefficients in an arbitrary module $M$, is analyzed in \cite{Degrijs-Petrosyan} for various $\frk$,  but unlike here, $\frk $ is required to act trivially on $M$.}.

\section{Filtration of the cochain complex}\label{improve}
Suppose for the moment that $\h$ is an arbitrary subalgebra of $\g$  and $M$ a $\g$-module.  If $F^n \in C^n(\g,M)$ vanishes whenever $r$ or more of its arguments lie in $\h$, then it is easy to check that  $\delta F^n$ vanishes whenever $r+1$ or more of its arguments lie in $\h$.
Let $\calF^jC^n(\g,M)$ be the space of those cochains which vanish whenever $n+1-j$ or more arguments lie in $\h$.  Then $\calF^jC^*(\g,M)$ is a subcomplex of $C^*(\g,M$), and we have the Hochschild-Serre descending filtration, \cite{Hoch-Serre:Lie},
\begin{equation}\label{HSfiltration}
C^*(\g,M) \; =\;\calF^0C^*(\g,M)\; \supset\; \calF^1C^*(\g,M)\;\supset\; \cdots \;\supset \calF^jC^*(\g,M)\;\supset \cdots.
\end{equation}
(For any fixed dimension $n$ the filtration terminates since $\calF^{n+1}C^n(\g,M)=0$.) This admits a stronger form when
$\h$ is a toral subalgebra of $\g$ acting torally on $M$.  Denote the homogeneous part of weight $\mbfz$ of $C^n(\g,M)$ by $C^n(\g,M)_{\mbfz}$ and the set of those  $F^n \in C^n(\g,M)_{\mbfz}$ which vanish whenever $r$ or more of its arguments lie in $\h$ by
$\calF_rC^n(\g,M)_{\mbfz}$. The cochain subcomplex $C^*(\g,M)_{\mbfz}$  inherits the filtration (\ref{HSfiltration}), but iteration of \eqref{cobound} shows, more strongly, that $\delta{\cal F}_rC^n(\g,M)_{\mbfz} \subset {\calF_r}C^{n+1}(\g,M)_{\mbfz}$.  
Also, if $F$ is a homogeneous $n$-cocycle of weight $\mbfz$, and  if we are given $h_1, \dots, h_r \in \h$, then $\iota_{h_r}\iota_{\h_{r-1}}\cdots\iota_{h_1}F$ is an $n-r$-cocycle whose cohomology class depends only on that of $F$ and the $h_i$. It is an alternating function of the latter, so we have a linear map
\begin{equation}\label{cochaindown}
\bigwedge\!\!{}^r\,\h\,\bigotimes C^n(\g, M)_{\mbfz} \longrightarrow C^{n-r}(\g,M)_{\mbfz},\quad r\le n.
\end{equation}
\noindent To put this in a form more usable later, observe that whenever $U, V,$ and $W$ are $\bfk$-spaces with $\dim_{\bfk}U < \infty$,  there is a canonical isomorphism
\begin{equation}\label{dual}
 \Hom(U\bigotimes V, W)\; \;\cong \;\; \Hom(V, U^{\vee}\bigotimes W),
 \end{equation}
where $U^{\vee}$ is the dual vector space to $U$. For suppose that $\dim U = q$, choose a basis $u_1,\dots,u_q$ and let $u_1^{\vee},\dots,u_q^{\vee}$ be the dual basis.  The isomorphism (\ref{dual}) sends a morphism $\phi:U\otimes V \to W$ to the morphism $\psi:V \to U^{\vee}\otimes W$ defined by sending  $v\in V$ to $ \sum_{i=1}^q u_i^{\vee} \otimes\phi(u_i \otimes v)$;
 this does not depend on the choice of basis.  For the inverse, suppose that $\psi:V\to U^{\vee}\otimes W$ is given.  If $v\in V$ then $\psi(v)$ is of the form $\sum_{i=1}^qu_i^{\vee}\otimes w_i$ and $\phi$ is defined by setting $\phi(u\otimes v) = \sum\langle u_i,v\rangle w_i$.

Suppose (as in $\mathfrak{sl}_N$) that $\dim \h = N\!-\!1$. Choose a basis  $\eta_1,\dots,\eta_{N\!-1}$, let $\eta_i^{\vee}, i = 1, \dots, N\!-1$ be the dual basis, and fix an integer $r$ with $0\le r\le N-1$. Then the $\eta_{i_1}\wedge\cdots\wedge \eta_{i_r}$ with $1\le i_1 < \cdots <i_r \le N\!-1$ form a basis for $\bigwedge^r \h$, where for $r=0$ this will mean $1 \in \bfk$. 
For simplicity, denote $\eta_{i_1}\wedge\cdots\wedge \eta_{i_r}$by $\eta_I$, where $I$ denotes the linearly ordered set of indices $i_1,\dots,i_r$. Then the $\eta_I^{\vee}=  \eta_{i_1}^{\vee}\wedge\cdots\wedge \eta_{i_r}^{\vee}$ form the dual basis for $\bigwedge^r\h^{\vee}$. If $F^n\in C^n(\g,M)_{\mbfz}$ then we will write  $\iota_{\eta_I} F^n$ for $\iota_{\eta_{i_r}}\iota_{\eta_{i_{r-1}}}\cdots\iota_{\eta_{i_1}}F^n$.  Applying \eqref{dual} now to \eqref{cochaindown} 
gives an epimorphism 
\begin{equation}\label{sigma}
\sigma:C^n(\g, M)_{\mbfz} \;\longrightarrow \;\bigwedge\!\!{}^r\,\h^{\vee} \otimes C^{n-r}(\g, M)_{\mbfz}
\end{equation}
sending $F^n \in C^n(\g,M)_{\mbfz}$ to $\sum_I\eta_I^{\vee}\otimes \iota_{\eta_I}F^n$. To see that it is onto, extend $\{\eta_1, \dots, \eta_{N-1}\}$ to an ordered basis of $\g$. Every $\eta_I^{\vee} \otimes f^{n-r}$ with $f^{n-r} \in C^{n-r}(\g, M)_{\mbfz}$ then has a preimage $F^n \in C^n(\g,M)_{\mbfz}$ defined by setting $F(g_1, \dots, g_n) = 0$
 when the set of  basis elements $\{g_1,\dots, g_n\}$ does not contain $\{\eta_{i_1}, \dots ,\eta_{i_r}\}$, and by $F(\eta_{i_1}, \dots ,\eta_{i_r}, g_{r+1},\dots,g_n) = 1$  when all the arguments are distinct basis elements and  are in the prescribed order. The kernel of $\sigma$ is $\calF_rC^n(\g,M)_{\mbfz}$, for it consists of those $F$ such that $\iota_{\eta_I}F = 0$ for all $I$ with $\#I = r$, i.e., those $F$ which vanish when $r$ or more of its arguments lie in $\h$. Fixing $r$ and considering all $n$, we have an exact sequence of complexes
\begin{equation}\label{dualexact}
0\;\longrightarrow \;\calF_rC^*(\g,M)_{\mbfz}\;\longrightarrow \;C^*(\g, M)_{\mbfz} \;\longrightarrow\; \bigwedge\!\!{}^r\,\h^{\vee} \otimes C^{*-r}(\g, M)_{\mbfz} \;\longrightarrow \;0.
\end{equation}
The coboundary operator of the quotient,  $\bigwedge\!\!{}^r\,\h^{\vee} \otimes C^{*-r}(\g, M)_{\mbfz}$ is, up to sign, just the coboundary operator of the second tensor factor.  For as $\bigwedge\!\!{}^*\h^{\vee}$ is identical with $C^*(\h, \bfk)$, the product   $\bigwedge\!\!{}^*\,\h^{\vee} \otimes C^{*-r}(\g, M)_{\mbfz}$ is a tensor product of two complexes, hence a complex, but $\h$ is Abelian, so the coboundary operator in the first tensor factor is zero. The short exact sequence \eqref{dualexact} gives rise, for every $r$, to a long exact sequence of cohomology groups, but we do not need this here.

\section{Toral semidirect products} \label{semidirect products}
If $\g = \h \ltimes \frk$ is a toral semidirect product, and if
$\h$ also acts torally on the $\g$-module $M$ then, as mentioned, the subcomplex  of $C^n(\g,M)_{\mbfz}$ of $C^*(\g,M)$ consisting of cochains of  weight $\mbfz$ relative to the weighting induced by $\h$ inherits the filtration (\ref{HSfiltration}). However, the filtration now actually arises from a gradation on $C^*(\g, M)_{\mbfz}$, causing the associated Hochschild-Serre spectral sequence to collapse.

As a vector space, $\g$ is just the direct sum of $\h$ and $\frk$, so we may obtain an ordered basis of  $\g$ by first taking the ordered basis $\eta_1,\dots, \eta_{N-1}$ of $\h$, followed by any ordered basis $\kappa_1, \dots, \kappa_K$ of $\frk$.  A cochain $F \in C^n(\g,M)$ will be completely determined by its values when its arguments $g_1,\dots,g_n$ are taken from this basis of $\g$ and are in the prescribed order, something which we will henceforth always assume.

Following Hochschild-Serre \cite{Hoch-Serre:Lie},  recall that two $\g$-modules $M$ and $M'$ are said to be \emph{paired} to a third, $P$, if there is a bilinear map  $M\otimes M \to P$ such that, denoting the image of $m\otimes m$ by $m \smile m'$,  one has \,$[g, m\smile m'] =[g,m] \smile m + m\smile [g,m']$ \,for all $g \in \g, m\in M, m' \in M'$. If $F^r\in C^r(\g,M), G^s\in C^s(\g, M')$ then define $F^r\smile G^s \in C^{r+s}(\g, P)$ as follows. Let $I = (i_1,\dots, i_r)$ be an ordered subset of $\{1,\dots,r+s \}$,  set $g_I = g_{i_1}\wedge\cdots \wedge g_{i_r}$, and if $J$ is its complement, define $g_J$ similarly. Then $I \bigsqcup J$ is a permutation of $\{1,\dots, r+s\}$. Letting $\nu(I)$ denote its signum,  set
$$F^r\smile G^s(g_1,\dots,g_{r+s}) \;=\; \sum \nu(I)F^r(g_I)\smile G^s(g_J),$$ where the sum is over all partitions of $\{1,\dots,r+s\}$ into a disjoint union $ I\bigsqcup J$ with $\# I =r,\, \# J = s$.  Then 
$$\delta(F^r\smile G^s)\; = \;\delta F^r \smile G^s + (-1)^r F^r \smile \delta G^s,$$
from which  it follows that the cup product descends to cohomology.
Since $M$ and $M'$ are always paired to $M \bigotimes M'$ and \, $\bfk\bigotimes M = M$,  the coefficient ring $\bfk$, considered as a trivial $\g$-module, is always paired with any $\g$-module $M$ to the same $M$. The cup product defines  morphisms
\begin{equation}\label{cochainpairing0}
C^r(\g,\bfk)\otimes C^s(\g,M) \;\longrightarrow \;C^{r+s}(\g,M)
\end{equation}
in which cup products of cochcain of weight zero again have weight zero.
The morphism $\g \to \h$ induces a morphism 
$$\bigwedge\!\!{}^r\h^{\vee} \otimes C^{n-r}(\g, M) \longrightarrow \bigwedge\!\!{}^r\g^{\vee} \otimes C^{n-r}(\g, M),$$ 
in which cochains of weight zero are carried to cochains of weight zero. Combining this with \eqref{cochainpairing0}, gives  a morphism
$$ \rho:\bigwedge\!\!{}^r\,\h^{\vee} \bigotimes C^{n-r}(\g, M)_{\mbfz} \;\longrightarrow \;C^n(\g,M)_{\mbfz}.$$
\begin{theorem}\label{collapse}
The composite morphism
$$\bigwedge\!\!{}^r\,\h^{\vee} \bigotimes C^{n-r}(\g, M)_{\mbfz} \;\stackrel{\rho}\longrightarrow \;C^n(\g,M)_{\mbfz} \;\stackrel{\sigma}\longrightarrow\; \bigwedge\!\!{}^r\,\h^{\vee} \bigotimes C^{n-r}(\g, M)_{\mbfz}$$
is the identity; the sequence (\ref{dualexact}) splits.
\end{theorem}

\noindent \textsc{Proof.}  Since $\bigwedge\!\!{}^r\,\h^{\vee} \bigotimes C^{n-r}(\g, M)_{\mbfz}$ is spanned by elements of the form $\eta_I \otimes G$ with $G\in C^{n-r}(\g, M)_{\mbfz}$ and $\eta_I$ of the form $\eta_{i_1}\wedge\cdots\wedge \eta_{i_r}$, it is  sufficient to prove that $\sigma\rho$ is the identity on such an element. Set $\rho(\eta_I \otimes G) = F$.   With the preceding convention,  $F(g_1,\dots , g_n)$ vanishes unless  $g_i = \eta_{i_ i}, \dots, g_r= \eta_{i_r}$, in which case its value is $G(g_{r+1},\dots,g_n)$. Now recall that $\sigma F = \sum \eta_J^{\vee}\otimes \iota_{\eta_J}F$.  By the definition of $F$ we have  $\iota_{\eta_J}F= 0$ unless $J=I$ and $ \iota_{\eta_I}F = G$. $\Box$
\medskip

Let $C^{n-r:r}(\g,M)_{\mbfz}$ now denote the subspace of $C^n(\g,M)_{\mbfz}$ spanned by those those $n$-cochains which, when its arguments are chosen as above, vanish unless exactly $r$ of its arguments are amongst the $\eta_i$. Then $\delta C^{n-r:r}(\g,M)_{\mbfz} \subset C^{n-r+1:r}(\g,M)_{\mbfz}$. For if  $F = F^{n-r:r}\in C^{n-r:r}(\g,M)_{\mbfz}$ and its arguments are chosen from the basis elements, then $\delta F$ vanishes if fewer than $r$ are amongst the $\eta_i$.  On the other hand, if $h_1,\dots,h_{r+1}$ are amongst the $\eta_i$, then as $F$ is of weight $\mbfz$ one has $\iota_{h_1}\cdots \iota_{h_{r+1}}\delta F = (-1)^{r+1}\delta(\iota_{h_1}\cdots \iota_{h_{r+1}}F)$ from (\ref{cobound}), but the right side vanishes by hypothesis.  Since every $F\in C^n(\g,M)_{\mbfz}$ can be written uniquely as a sum of components in the various $C^{n-r:r}(\g,M)_{\mbfz}$  we have the following decomposition into a direct sum of subcomplexes.

\begin{theorem}
$$C^*(\g,M)_{\mbfz}\; = \;\bigoplus_{r=0}^{N-1}C^{*-r:r}(\g,M)_{\mbfz} \quad\Box$$
\end{theorem}

Recall that $\eta_I = \eta_{i_1}\wedge\cdots \wedge \eta_{i_r}$ for an $r$-tuple of integers $1\le i_1 < i_2 \cdots < i_r \le N\!-\!1$. With $n$ fixed, similarly define $\kappa_J$ for a set of $n-r$ distinct integers between 1 and $K$.  When $F\in C^{n-r:r}$ and  arguments are restricted to the chosen  basis elements, one can express its value simply as $F(\eta_I,\kappa_J)$ for suitable $I$ and $J$. If also $F\in C^{n-r:r}_{\mbfz}$ then $\iota_{\eta_I}F = \iota_{i_r}\iota_{i_{r-1}}\cdots\iota_{i_1}F$ is an $n-r$ cochain of weight $\mbfz$ which vanishes when any argument is in $\h$, and so may be viewed as an element of $C^{n-r}(\frk,M)_{\mbfz}$.  We therefore have a cochain morphism (up to sign)
$$\bigwedge\!\!{}^r\h \bigotimes C^{*-r:r}(\g,M)_{\mbfz}\; \longrightarrow \;C^{*-r}(\frk,M)_{\mbfz},$$
where on the left the coboundary operator operates only on $C^{*-r:r}(\g,M)_{\mbfz}$.
From the  preceding section, this may be identified with a cochain morphism
\begin{equation*}
\phi: C^{*-r:r}(\g,M)_{\mbfz}\; \longrightarrow\; \bigwedge\!\!{}^r\h^{\vee} \bigotimes  C^{*-r}(\frk,M)_{\mbfz},
\end{equation*}
where now on the right the coboundary operator operates only on $C^{*-r}(\frk,M)_{\mbfz}$.
However, we also have a cochain morphism
\begin{equation*}
\psi:\bigwedge\!\!{}^r\h^{\vee} \bigotimes  C^{*-r}(\frk,M)_{\mbfz}\;\longrightarrow  \;C^{*-r:r}(\g,M)_{\mbfz}
\end{equation*}
defined as follows. If $\xi \in \bigwedge\!\!{}^r\h^{\vee}, f \in C^{n-r}(\frk,M)_{\mbfz}$, then to define $F = \psi(\xi\otimes f) \in C^{n-r:r}(\g, M)_{\mbfz}$ we only have to give its values when its arguments are amongst the chosen basis element of $\h$ and $\frk$.  Let the set of arguments be written, as above, in the form $\eta_{I},\kappa_{J}$, where $I$ is an $r'$-tuple of integers and $J$ an $s'$-tuple, with $r'+s' = \dim \mathfrak{g}$. Then set $F(\eta_{I},\kappa_{J})=0$ unless  $r'=r$, in which case set  $F(\eta_I,\kappa_J) = \langle \xi, \eta_I\rangle f(\kappa_J)$. This $\psi$ is the inverse of $\phi$. With $n$ fixed, summing over $r$ gives
$$C^n(\g,M)_{\mbfz}\; =\; \bigoplus_{r=0}^n \bigwedge\!\!{}^r \h^{\vee}\bigotimes C^{n-r}(\frk,M)_{\mbfz},$$
so there is an isomorphism of complexes
\begin{equation*}
C^*(\g,M)_{\mbfz}\; \cong \; \bigwedge\!\!{}\h^{\vee} \bigotimes  C^*(\frk,M)_{\mbfz}.
\end{equation*}

\noindent Taking cohomology yields the following.
\begin{theorem} \label{semidirect}
Let $\g= \h\ltimes\frk$ be  a semidirect product where $\h$ is Abelian and acts torally both on $\g$ and on a $\g$-module $M$. Then with respect to the weighting induced by the action of $\h$ we have
$$H^*(\g,M)\; \cong \; \bigwedge\!\!{}\h^{\vee} \bigotimes  H^*(\frk,M)_{\mbfz} \; =\; H^*(\h,\bfk)\bigotimes  H^*(\frk,M)_{\mbfz}. \quad \Box $$
\end{theorem}
Note that $H^0(\frk,M)_{\mbfz}$ consists of the invariants of $M$ under the operation of all of $\g$, for it consists of the invariants of $M$ under $\frk$ which are also invariant under $\h$, but as a vector space, $\g$ is the direct sum of $\frk$ and $\h$

\section{Lie poset subalgebras of $\mathfrak{sl}(N)$}
Let $\b$ again be the Borel subalgebra of  $\mathfrak{sl}(N)$ and $\n$ now be the ideal of $\b$ consisting of all strictly upper triangular matrices.  While $\h$ does operate to give a decomposition of the cohomology,  as $\h$ is not contained in $\n$, the cohomology does not necessarily reduce to the weight $\mbfz$ part \footnote{The cohomology of $\n$ with coefficients in an arbitrary module has been computed by Kostant, \cite{Kostant:Borel}, but the nature of the decomposition does not seem to have been considered by him.}. Recall  that $\eta_i = e_{ii}-e_{i+1,i+1}$. 
\begin{theorem} \label{Lie poset}
Let $\g$ be  a subalgebra of $\mathfrak{sl}(N)$ with $\b \supset \g \supset \h$ and set $\frk = \g \cap \n$. Then \,(i) $\g = \h \ltimes \frk$, and (ii) \,$\frk$ is spanned by those $e_{ij}$ which it contains; defining a partial order $\mathcal P$ on $\{1,\dots,N\}$ by setting $i\preceq j$ if either $i = j$ or $e_{ij} \in \frk$, \, (iii)\, $\g$ is the Lie poset algebra $\g(\mathcal P)$.
\end{theorem}

\noindent \textsc{Proof.} The first assertion is immediate from the fact that an upper triangular matrix is uniquely a sum of an element of $\h$ and an element of $\n$. For the second it is sufficient to show that if some linear combination $a=\sum_{k=1}^r c_ke_{i_k,j_k}, i_k < j_k, c_k \ne 0$ is in $\frk$ then at least one of the $e_{i_k,j_k}$ is already in $\frk$; this will imply that all are in $\frk$. If not, suppose that the given  $a \in \frk$ is one with minimal $r$ having no summand $e_{i_k,j_k}$ in $\frk$; surely $r\ge 2$.  Then $[\eta_{i_1}, a]$ is not a multiple of $a$ but is a linear combination of the same summands $e_{i_k,j_k}, k = 1,\dots,r$, so there is a linear combination of $a$ and $[\eta_{i_1},a]$ which is not zero and contains no more than $r-1$ of these summands. One of them is consequently already in $\frk$, a contradiction. The last assertion follows. $\Box$
 \medskip
 
The description of Lie poset algebras in Theorem \ref{Lie poset}  is meaningful for all Lie algebras of Chevalley type. This suggests that the results that follow may carry over in some way to such algebras.
\medskip

From this point on we assume that the characteristic of \,$\bfk$ is greater than $N$.
\medskip

Since the $e_{ij}$ in $\mathfrak{sl}(N)$ are all simultaneous eigenvectors for the operations of $\h$ they determine elements of $\h^{\vee}$. The weight defined by $e_{ij}$ will be denoted $\w_{ij}$, so $[h,e_{ij}] =\w_{ij}(h)e_{ij}$ for $h\in\h$. In what follows we use the fact that every $\w_{ij}$ with $i<j$ is a sum of weights of simple positive roots: $\w_{ij} = \w_{i, i+1} + \w_{i+1,i+2} + \cdots + \w_{j-1,j}$.
\begin{theorem}\label{uniqueness lemma}
If $e_{i_1,j_1},\dots,e_{i_k,j_k},\, i_r < j_r, r = 1, \dots, k$ are distinct elements of $\mathfrak{sl}(N,\bfk)$, then (i)\,
$\sum_{r=1}^k \w_{i_r,j_r}\ne \mathbf{0}$ and  (ii)
$\sum_{r=1}^k \w_{i_r,j_r}\ne \w_{ij}$ for any $e_{i,j}$ unless the $e_{i_r,j_r}$ can be so ordered that $i = i_1, j_r = i_{r+1}, r= 1, \dots, k-1$ and $j_k = j$.
\end{theorem}

\noindent \textsc{Proof.} 
Suppose (i) were false. Reordering if necessary, we may assume that $i_1$ is minimal amongst the $i_r$. Then we may assume that $i_1 = 1$, else we could reduce the value of $N$.  
Consider the summands on the right of the form $\w_{1j}$. Since $[\eta_1, e_{12}] = 2$ while $[\eta_1, e_{1j}] = 1$ for $j >1$, writing all $\w_{i_r,j_r}$ as sums of weights of simple positive roots, one sees that the sum of these can not vanish. For $\w_{12}$ appears, but can not appear with coefficient greater than $N$ since there are no more than $N-1$ distinct possible summands of the form $\w_{1j}$; as we have assumed that $p > N$ it can not be zero modulo $p$. 

For (ii), suppose that we have an equality of the kind given. If $i_1$ is minimal amongst the $i_r$ that appear on the left then $i_1 = i$, for it clearly can not be greater, but can not be less since, by the same argument as before, the summands of the form $\w_{i_1,j_r}$ could not cancel.  Similarly there can not be more than one summand of the form $\w_{i,j_r}$, so $i_1 = i$,  implying that $r = 1$. If now $i_2 > i_1$ is minimal amongst the remaining $i_r$ then we must have $i_2 \ge j_1$.  Otherwise, writing every weight on the left as a sum of weights of simple positive roots, observe that $\w_{i_2,i_2+1}$ would occur at least twice, but its multiplicity can not then be 1 modulo $p$ since there can be no more than $p-1$ elements of the form $\w_{i_2, j_r}$ on the left. Continuing, we see that after possible reordering we must have $i = i_1, j_1 \le i_2, j_2 \le i_3, \dots$, from which it is clear that the assertion must hold.$\,\Box$\medskip

It follows that $H^n(\frk,\bfk)_{\mbfz} = 0$ for $n>0$, for a cochain $F$ of weight $\mbfz$  evaluated on elements of $\frk$ must have non-zero weight, but all elements of the coefficient module have weight $\mbfz$.  On the other hand, $H^0(\frk,\bfk)_{\mbfz} = \bfk$, so with Theorem \ref{semidirect} one can reproduce the result of \cite{G:xLie}.
\begin{theorem}\label{trivial module}  
\qquad $H^*(\g(\calP), \bfk) \;\cong\; \bigwedge^*\h^{\vee} \; = \; H^*(\h, \bfk).  \qquad \Box$
\end{theorem}
Theorem \ref{uniqueness lemma} will permit us to identify $H^n(\frk,\g(\calP))_{\mbfz}$ with $H^n(\Sigma, \bfk)$ for $n>0$, in turn allowing the application of Theorem \ref{semidirect}.  A non-degenerate $n$-simplex of the partially ordered set  $\{1,\dots,N\}$ can be identified with an ordered $n+1$ tuple $(i_0,i_1,\dots,i_n)$ of integers with $i_0\prec i_1 \prec \cdots \prec i_n$ in the partial order induced by $\g$, and an $n$-cochain in $C^n(\Sigma, \bfk)$ can then be considered as a function $f^n(i_0,i_1,\dots,i_n)$ from such $n$-tuples to $\bfk$, where $\Sigma$ is the associated simplicial complex. Define a mapping $\Phi^n: C^n(\Sigma, \bfk) \to C^n(\frk,\g(\calP))_{\mbfz}$ for all $n>0$. A cochain is uniquely determined by its values when all its arguments are basis elements. If $f^n \in C^n(\Sigma, \bfk)$ then set the value of $\Phi^n f^n$ equal to zero if any argument is an $\eta_i$; it then vanishes if any argument is in $\h$.
If $i_0,\dots,i_n \in \{1,\dots, N\}$  with $i_0\prec i_1 \prec \cdots \prec i_n$ (so by hypothesis $e_{i_0,i_1},e_{i_1,i_2},\dots,e_{i_{n-1},i_n}\in \frk$), then set  $(\Phi f^n)(e_{i_0,i_1},e_{i_1,i_2},\dots,e_{i_{n-1},i_n}) = f^n(i_0,\dots,i_n)e_{i_0,i_n}$ and extend the definition so that $\Phi f^n$ is an
alternating multilinear function of these arguments. Finally,  set $(\Phi f^n)(e_{i_0,j_0},e_{i_1,j_1},\dots,e_{i_{n-1},j_{n-1}}) = 0$ if the arguments can not be so reordered that $j_0 = i_1, j_1 = i_2, \dots, j_{n-2} = i_{n-1}$.    Then $\Phi f^n$ is homogeneous of weight $\mathbf{0}$ relative to the weighting induced by the toral subalgebra $\h$ of $\g(\calP)$, hence an element of $C^n(\frk, \frk)_{\mbfz}$, which is identical with 
$C^n(\frk,\g(\calP))_{\mbfz}$ for $n>0$ . Theorem \ref{uniqueness lemma} asserts that $\Phi^n$ is onto.

When $n=0$ we must make a minor modification (often introduced to avoid exceptional cases, cf, e.g., in Alexander duality). Augmenting $\Sigma$ to $\Sigma^+$ by adjoining  a single simplex of dimension $-1$ to serve as the boundary of every $0$-simplex one obtains the reduced cohomology $H^*(\Sigma^+, \bfk)$. When char $\bfk = 0$ this coincides with $H^*(\Sigma,\bfk)$ in every dimension except 0. While this may fail for positive characteristic, the following shows that it continues to hold for the cases we are considering.
\begin{theorem} With the preceding notation,
$H^0(\Sigma^+, \bfk)$ can be naturally identified with the center $\ctr$ of $\g$, while
for $n \ne 0 $ one has $H^n(\Sigma^+, \bfk)\; = \;H^n(\Sigma, \bfk).\quad $ 
\end{theorem}

\noindent \textsc{Proof.} The center of  $\g$ clearly is contained in $\h$. A diagonal matrix may be viewed as a function on $1,\dots,N$ viewed as $0$-simplices. It  is in $\ctr$ precisely when it is constant on their homology classes, so these elements of $\h$ are just the functions on the homology classes, hence constitute the cohomology in dimension zero, proving the first assertion.  Since $C^*(\Sigma^+, \bfk)$ and $C^*(\Sigma, \bfk)$ differ only at dimension 0,
 we certainly have $H^n(\Sigma^+, \bfk) = H^n(\Sigma, \bfk)$ for $n>1$. To show that also $H^1(\Sigma^+, \bfk) = H^1(\Sigma, \bfk)$ we must show that $\delta C^0(\Sigma^+, \bfk) = \delta C^0(\Sigma, \bfk)$.  However, every 0-cochain $f:\{1,\dots, N\} \to \bfk$ in $C^0(\Sigma, \bfk)$ can be written uniquely as a sum $f= f' + f''$ with $f' \in C^0(\Sigma^+, \bfk)$ and $f''$ having the same constant value on all $i \in \{1,\dots,N\}$: set $f'(i) = f(i) - (1/N)\sum_{j=1}^N f(j)$ for all $i$ and $f''(i) = (1/N)\sum_{j=1}^N f(j)$ all, $i$. This is well-defined since the characteristic is greater than $N$. Then $\delta f'' = 0$,  proving the second assertion. Since the unique simplex of dimension $-1$ of $\Sigma^+$ is always a coboundary it does not contribute to the cohomology.  $\Box$
 \smallskip

\begin{theorem}\label{Phi} The mapping \, $\Phi: C^*(\Sigma^+, \bfk) \to C^*(\frk,\g(\calP))_{\mathbf 0}$ is a cochain isomorphism.
\end{theorem}

\noindent \textsc{Proof.} We must show that $\delta \Phi f^n= \Phi \delta f^n$ for all $f^n \in C^n(\Sigma^+, \bfk)$. For $n > 0$ it is sufficient to show that the two sides coincide when the  arguments 
 are amongst the $e_{ij}, i \prec j$.  Since $\delta \Phi f^n$ is a homogeneous cocycle of weight $\mathbf{0}$ it must vanish, by Theorem \ref{uniqueness lemma}, unless the arguments can be reordered to be of the form $(e_{i_0,i_1}, e_{i_1,i_2},\dots, e_{i_{n-1},i_n})$ with $i_0\prec i_1 \prec \cdots \prec i_n$, in which case the value must be a multiple of $e_{i_0,i_n}$. It follows, cf. \!Theorem \ref{uniqueness lemma}, that the only non-zero terms in $(\delta \Phi f^n)(e_{i_0,i_1}, e_{i_1,i_2},\dots, e_{i_{n-1},i_n})$ can be ones of the form
\begin{multline*}
[e_{i_0,i_1}, (\Phi f^n)(e_{i_1,i_2},\dots, e_{i_{n-1},i_n})],\,[e_{i_{n-1},i_n},(\Phi f^n)(e_{i_0,i_1},\dots, e_{i_{n-2},i_{n-1}}]\quad \text{ and }\\
(\Phi f^n)([e_{i_{r-1},i_r},e_{i_r,i_{r+1}}], e_{i_0},\cdots,\hat e_{i_{r-1},i_r}, \hat e_{i_r,i_{r+1}}, \dots e_{i_{n-1},i_n}),\quad r=1,\dots,n-1.
\end{multline*} 
It only remains to consider the signs with which these terms appear; examining them shows that indeed
 \begin{multline*} 
 \delta (\Phi f^n)
(e_{i_0,i_1}, e_{i_1,i_2},\dots, e_{i_{n-1},i_n})\quad = \quad((\delta f^n)(i_0,\dots,i_n))e_{i_0,i_n} \\= \Phi(\delta f^n)(e_{i_0,i_1}, e_{i_1,i_2},\dots, e_{i_{n-1},i_n}).\qquad
\end{multline*}
For the case $n=0$, recall ($\S \ref{semidirect products}$) that $H^0(\frk, \g(\calP))_{\mathbf 0}$ is  just the center $\ctr =\ctr(\g(\calP))$ of $\g(\calP)$. This is spanned by those diagonal matrices $h$ such that $h(i) = h(j)$ whenever $e_{ij} \in \frk$.
After augmenting $\Sigma$, those 0-cochains having the same constant value on all $i \in \{1,\dots,N\}$ are now coboundaries. A 0-cocycle of $\Sigma^+$ is a cochain $f$ such that $f(i) = f(j)$ whenever $e_{ij} \in \frk$. For each there is now a unique cohomologous cochain with $\sum_{i=1}^Nf(i) = 0$, so we also now have an isomorphism, $\ctr(\mathfrak{g}(\calP)) \cong H^0(\Sigma^+,\bfk). \quad  \Box$ 
\medskip

This yields our final theorem.
\begin{theorem}\label{final} 
If $\cal{P}$ is a partial oder on $\{1,\dots, N\}$ and char $\bfk > N$ or \,$0$, then
$$\qquad H^*(\mathfrak{g}(\calP),\mathfrak{g}(\calP))\;= \;\bigwedge \h^{\vee} \bigotimes H^*(\Sigma^+, \bfk)\; =\; H^*(\h, \bfk) \bigotimes H^*(\Sigma^+, \bfk)\qquad \Box$$
\end{theorem}

\section{Some examples of deformations}

With Theorem \ref{final} we can return to the discussion in the Introduction.  We have finally proven that if $\cal P$ is a finite poset, then writing  $\g$ for  $\g(\cal{P})$ and $\Sigma$ for $\Sigma(\cal P)$, the space $H^2(\g, \g)$ of infinitesimal deformations of $\g$  is
\begin{equation}\label{inf_defs}
H^2(\g,\g) \;=\; (\bigwedge{}\!^2\,\h^{\vee} \bigotimes \ctr)\; \;\bigoplus\; \;(\h^{\vee} \bigotimes H^1(\Sigma,\bfk))\; \;\bigoplus\; \;(H^2(\Sigma, \bfk)).
\end{equation}
The infinitesimals in the first summand (``type (2,0)'') involve only elements of $\h$ since any central element of a Lie poset algebra $\g$ must lie in $\h$.  Those infinitesimals in the last summand (``type (0,2)'') involve only elements of $\frk$. The elements in the first and last summands do not change the operation of $\h$ on $\frk$; the infinitesimals in the middle summand (``type (1,1)'') do. Some simple examples follow of deformations whose infinitesimals are of each of these three kinds. 

Observe from \eqref{inf_defs} that the necessary and sufficient condition for a Lie poset algebra to be absolutely rigid, i.e., to have no infinitesimal deformations,  is the simultaneously vanishing of $\ctr, \, H^1(\Sigma, \bfk)$, and $H^2(\Sigma, \bfk)$. An associative poset algebra $A(\mathcal P )$ may be absolutely rigid while its associated Lie poset algebra $\g(\mathcal P)$ may allow non-trivial deformations.  In the simplest case, let the partial order on $\{1,\dots,N \}$ be vacuous. Then $A(\mathcal P )$ is the algebra of all diagonal $N \times N$ matrices, a direct sum of $N$ copies of $\bfk$, hence separable. Having trivial cohomology, it is absolutely rigid.  However, 
$\h$ is Abelian and therefore deformable to any algebra of the same dimension.  By contrast, suppose that $\cal{P}$ is $\{1,\dots,N\}$ with partial order identical with the  linear order. Then $\ctr  =0$ and  $H^*(\Sigma^+, \bfk) = 0$, so just as in the associative case, $\g(\cal{P})$ is absolutely rigid. Since all Lie poset algebras are solvable this provides a simple example of a solvable Lie algebra which is absolutely rigid. It is not known, both in the associative and the Lie cases, whether a nilpotent algebra can be rigid, let alone absolutely rigid.

One can construct small posets ${\cal P}^n$ whose corresponding simplicial complexes have geometric realizations which are $n$-spheres, \cite[\S 15]{GS:Monster}. Start with the zero sphere, $S^0$, which consists of just two points, and repeatedly take two-point suspensions. For $S^0$ one has  ${\cal P}^1$ = $\{1,2 \}$ with vacuous partial order,  for $S^1$  take its two-point suspension ${\cal P}^2$ = $\{1,2,3,4\}$ with $1 \prec 3, 4; 2\prec 3,4$ and no other relations, for  $S^2$ take ${\cal P}^3$ = $\{1,2,3,4,5,6 \}$ with $1\prec 3,4; 2\prec 3,4 ; 3\prec 5, 6; 4 \prec 5,6$ and no relations other than those following from these, and so on. Denote the corresponding associative poset algebras by $A(S^n)$ and the Lie poset algebras by $\g(S^n)$.
Since  $H^2(A(S^n), A(S^n)) = H^2(S^n, \bfk)$, which vanishes  except when $n=2$, it follows that $A(S^n)$ is rigid for $n\ne 2$. When $n=2$ there is exactly one infinitesimal deformation, and as $H^3(A(S^2), A(S^2)) =0$, it is unobstructed. Therefore $A(S^2)$ has a one-parameter family of non-trivial deformations.  Assuming that the characteristic of  $\bfk$ is at least 6, Theorem \ref{final} asserts that $H^2(\g(S^2), \g(S^2))$ is likewise one-dimensional. Although $H^3(\g(S^2), \g(S^2))$ does not vanish, the existing single infinitesimal deformation of $\g(S^2)$ is not obstructed and in fact does give rise to a one parameter family of deformations. For while a  deformation of an associative algebra $A$ always gives rise to a  deformation  of its commutator Lie subalgebra which could be trivial, nevertheless, although $\g(S^2)$ is not just the commutator subalgebra of $A(S)$, the deformation is essentially the induced one and is non-trivial because its infinitesimal, which is of the type (0,2), is non-trivial. Theorem \ref{final} also implies that $\g(S^n)$ is absolutely rigid for $n> 2$.

Deformations induced by infinitesimals of the type (2,0) are necessarily jump deformations, since after  deformation $\h$ will no longer be commutative. The cohomology of the algebra will have changed and the original algebra can not be recovered as a deformation of the new one.  In the extreme case, where $\g$ is reduced to $\h$ and the ideal $\frk$ vanishes, $\g$ is Abelian and can deform to any Lie algebra of the same dimension, in particular, to another Lie poset algebra.  This shows that the underlying topology, i.e., that of the geometric realization of the nerve, which was originally discrete, can cease to be so.  A small example where it does not change, however, is given by the three-dimensional subalgebra $\g$ of $\mathfrak{sl}(3)$ spanned by $\eta_1= e_{11}-e_{22},\, \eta_2 = e_{22}-e_{33}$, and $e_{12}$. Its infinitesimal deformations are all of the type (2,0). The center is non-trivial and is spanned by $ \eta_1+2\eta_2 = e_{11}+e_{22}-2e_{33} $. Denoting this by $c$, the single element $(\eta_1 \wedge \eta_2) \otimes c$ spans $H^2(\g,\g)$ and gives rise to a deformed Lie product $[-,-]_*$ in which we have $[\eta_1,\eta_2]_* = tc$ with all other products of basis elements unchanged. The power series one normally encounters is here just a single term, but infinitesimals of the type (2,0) in general can be obstructed.  Here, setting $t=1$, the Cartan subalgebra $\h$ of $\g$ has been deformed to the unique non-abelian two-dimensional Lie algebra. Viewing $\h$ as the algebra of functions on the vertices of $\Sigma$, which here is just a one-simplex,  its multiplication has been deformed, becoming non-commutative, but the evaluation of an $h \in \h$ on any $e_{ij}$, i.e., the corresponding eigenvalue of $\operatorname{ad}h$, has not changed.

Infinitesimals of the type (1,1) are of the form $z = \sum \xi_i \otimes D_i$ where $\xi_i \in \h^{\vee}$ and $D_i$ are derivations of $\frk$ which commute with the operation of $\h$ and therefore send every $e_{ij} \in \frk$ to a multiple of itself. To first order, the deformation induced by $z$ is given by $[h, e_{ij}]_* = [h, e_{ij}] + t\sum\langle \xi_i,h\rangle e_{ij}$ for all $h \in \h$; products of two elements of $\h$, as well as of two elements of $\frk$, remain unchanged. Such infinitesimals may also meet obstructions except, as in the complex analytic case, when that part of the third cohomology group vanishes which contains the possible obstructions to the given kind of infinitesimal. In the complex case this means that  $H^2( \cal X, \cal T)$ = 0, e.g., when the manifold is a Riemann surface. (Note that ${\cal T} = \bigwedge{}\!^1{\cal T}$, contributing one dimension.) Here the condition is that $\h^{\vee} \bigotimes H^2(\Sigma, \bfk) = 0$, but as $\h^{\vee}$ is never zero, it becomes simply that $H^2(\Sigma, \bfk) = 0$. This is the case, for example, when ${\cal P} = {\cal P}^1$. 

As an illustration, we compute the deformation associated with one particular non-trivial cocycle. In $\g(S^1)$ the ideal $\mathfrak k$ is spanned by $e_{13}, e_{14},e_{23},e_{24}$, and is Abelian. As  $\Sigma(\mathcal P )$ has no 2-simplices, any 1-cochain $f$ there is already a 1-cocycle.  The condition for $f= f^1$ to be a coboundary is that $f(1,2)-f(2,3)+f(2,4) - f(1,4) = 0$; viewing $f$ as an alternating function of its arguments, this says that $f(1,2) +f(2,3) +f(3,4) + f(4,1) = 0.$ So, for example, choosing $f(1,2) = 1, f(1,3) = f(2,3) = f(2,4)= 0$ gives a 1-cocycle whose class generates $H^1(\Sigma(\mathcal P),\bfk)$.  The Cartan subalgebra $\h$ of $\g(S^1)$ is 3-dimensional with basis $\{\eta_i, i = 1,2,3\}$ in the notation of the preceding section. Let $\{\eta_i^{\vee}, i = 1,2,3\}$ be the dual basis.  Then $ z =\eta_1^{\vee}\otimes \Phi f$ is a non-trivial 2-cocycle of $\g(S^1)$. As basis for $\g(S^1)$ we  have  $\{\eta_1, \eta_2,\eta_3, e_{13},  e_{14}, e_{23}, e_{24}\}$.  The deformed product $[-,-]_*$ induced by $z$ on these basis elements is exactly the same as the original except in one instance: one has $[\eta_1,e_{12}]_* = (2+t)e_{12}$. This defines the entire deformation, but note, for example that we now have $[e_{11} -e_{33}, e_{12}]_* = (1+t)e_{12}$, since we must write $e_{11}-e_{33}$ as $\eta_1 + \eta_2$. We can specialize $t$ to values in $\bfk$ to obtain an algebra defined over $\bfk$. For almost all values one can recover the original algebra as a specialization of a deformation of the new one, but it is not clear if this is possible when, for example, we set $t$ equal to -2.  Similarly, when a complex variety is deformed, it may become singular or degenerate for some values of the deformation parameter.

To understand the mysterious  infinitesimals of the type (0,2), suppose that we have an associative poset algebra $A(\cal P)$, and an element of $H^2(\Sigma({\cal P}), \bfk)$ represented by a 2-cocycle $f$. When $i\prec j \prec k$, denote by $f_{ijk}$ the value of $f$ on the 2-simplex these define. Those $e_{ij}$ in $\frk$, together with 0, form a semigroup; one may view $f$ as a 2-cocycle of this semigroup with  values in the additive group of $\bfk$  by setting $f(e_{ij}, e_{j'k}) =  f_{ijk}$ if $j = j'$ and 0 otherwise. If $\bfk$ is $\R$ or $\C$, then we can define a one-parameter family of deformations of the ideal $\frk$ of $A$ by setting $e_{ij}\star e_{jk} = (\exp t\!f_{ijk}) e_{jk}$, extended linearly. This is associative, and can be extended to all of $A$ by leaving unchanged products of any two elements in which one is an $e_{ii}$. Up to equivalence, the deformation so obtained depends only on the class of $f$. It induces a deformation of the associated Lie poset algebra $\g(\cal P)$, where the class of  $f$ is now a differential of the type (0,2).  Replacing $t$ by $i\!t$, the deformation becomes a `phase' periodic in $t$. 

The idea of exponentiating an ``additive'' 2-cocycle of a semigroup to get a ``multiplicative'' one, thereby yielding a family of deformations, may have first appeared in \cite{G:Path}, cf. also \cite{G:deBroglie}. A change in structure defined by an arbitrary  multiplicative 2-cocycle need not be a deformation in the technical sense, but the exponential of an additive 2-cocycle is. It can not be a jump deformation and one can recover the original algebra as a specialization of a deformation of the new one, in sense undoing the deformation.  
When $\bfk$ has finite characteristic, one can construct a formal family by use of the Artin-Hasse exponential, but it may not be possible then to specialize the deformation parameter $t$ to an element of the ground field $\bfk$.

As an illustration, observe that the only non-trivial infinitesimal deformations of $\g = \g(S^2)$ are of the type (0,2), stemming from the second cohomology of the 2-sphere $S^2$. When $\bfk$ is either $\R$ or $\C$,  the exponential of a non-trivial 2-cocycle $f$  of $\Sigma({\cal P}^2)$ yields a non-trivial deformation of $\g$. In this example, however, it is unnecessary to require that $\bfk$ be either $\R$ or $\C$, since $H^3(\Sigma, \bfk) = 0$. As all products of three of the $e_{ij}$ in $\g(S^2)$ vanish, every 2-cochain of $\Sigma({\cal P}^2)$ is a 2-cocycle, and in place of the full exponential we can simply use $F_{ijk}(t) = 1 + tf_{ijk}$  as a family of multiplicative 2-cocycles to obtain a one-parameter family of deformations; one sets $[e_{ij}, e_{jk}]_* = F_{ijk}(t)e_{ik}$, other products of basis elements being unchanged.  Since   $\dim H^2(\Sigma({\cal P}^2), \bfk) = 1$,  any 2-cochain $f$ which is not a coboundary will produce, up to equivalence, the unique one-parameter family of deformations which $\g(S^2)$ allows.

 In each of our examples, the deformation has produced an algebra which is no longer a Lie poset algebra, suggesting that there is some larger category  in which the theory should be set. This problem already arises in the case of associative poset algebras, as the last example shows.  Such algebras, however, are particular examples of \, ``Tic-Tac-Toe'' algebras in the sense of  Mitchell, \cite{Mitchell:Categories}, and any deformation of a Tic-Tac-Toe algebra is again a Tic-Tac-Toe algebra, \cite{GS:Triangular}, so the category of \,Tic-Tac-Toe algebras is stable under deformation. This raises the  question of whether there is a stable category of Lie algebras containing the category of Lie poset algebras as a full subcategory. A much more intriguing question, however, is whether the similarity between the deformation theories of complex analytic manifolds and Lie poset algebras actually stems from their being special cases of something more general.

 \nocite{HazGer}


\begin{thebibliography}{10}

\bibitem{Ch-E}
C. ~Chevalley and S.~Eilenberg.
\newblock{Cohomology theory of Lie groups and Lie algebras}.
\newblock{\em Trans. Amer. Math. Soc.}, 63:85--124, 1948.

\bibitem{Degrijs-Petrosyan}
D.~Degrijs and N. Petrosyan.
\newblock{On the cohomology of split Lie algebra extensions}.
\newblock{\em J. Lie Theory}, 22:1--15, 2012.

\bibitem{DitoSternheimer:Genesis}
G.~Dito and D.~Sternheimer.
\newblock{Deformation quantization: genesis, developments
and metamorphoses}, pp. 9--54 in
\newblock{G. Halbout, ed.,{\em Deformation Quantization, IRMA Lectures in Math. Theoret. Phys.
1},}
\newblock{Walter de Gruyter, Berlin} 2002.
\newblock{Cf. {\texttt arXiv:0201168v1 [math.QA] 18 Jan 2002},} 41 pp.

\bibitem{FN:Stability}
A.~Fr\"olicher and A.~Nijenhuis.
\newblock{A theorem on stability of complex structures}.
\newblock{\em Proc. Nat. Acad. Sci. (USA)}, 43(2): 239-241, 1957.

\bibitem{G:xLie}
M.~Gerstenhaber.
\newblock{A note on the cohomology of Lie algebras}.
\newblock{\texttt arXiv:1208.0350v1 [math.RT] 1 Aug 2012}, 4 pp.

\bibitem{G:Path}
M.~Gerstenhaber.
\newblock{On the deformation of path algebras}.
\newblock{\texttt arXiv:1306.4939v2[math-ph] 7 Jul 2013}, 11 pp.

\bibitem{G:deBroglie}
M.~Gerstenhaber.
\newblock{Path algebras and de Broglie waves}.
\newblock{\texttt arXiv:1403.3429v1 [math-ph] 13 Mar 2014}, 13 pp.

\bibitem{G:Self-dual}
M.~Gersenhaber.
\newblock{Self-dual and quasi self-dual algebras}
\newblock{\em Israel J. Math.} 200:193-211, 2014. 

\bibitem{GG:WeylCohom}
M.~Gerstenhaber and A.~Giaquinto.
\newblock{On the cohomology of the Weyl algebra, the quantum plane, and the $q$-Weyl algebra}.
\newblock{\texttt arXiv:1208.0346 [math.QA] 1 Aug 2012}, 16 pp.

\bibitem{GS:SC=HC}
M.~Gerstenhaber and S.~D.~Schack
\newblock{Simplicial cohomology is Hochschild cohomology}.
\newblock{\em J. Pure Appl. Algebra}, 30:143--156, 1983.

\bibitem{GS:Monster}
M.~Gerstenhaber and S.~D. Schack.
\newblock {Algebraic cohomology and deformation theory}, pp. 11--264
\newblock in [14].

\bibitem{GS:Triangular}
M.~Gerstenhaber and S.~D. Schack.
\newblock {Triangular algebras}, pp. 447--494
\newblock in [14].

\bibitem{Goldberg:Lie}
S.~I. Goldberg.
\newblock {On the Euler characteristic of a Lie algebra}.
\newblock {\em Amer. Math. Monthly}, 62:239--240, 1955.

\bibitem{HazGer}
M.~Hazewinkel and M.~Gerstenhaber, editors.
\newblock {\em {Deformation Theory of Algebras and Structures and
  Applications}}, volume 247 of {\em NATO ASI Science Series}.
\newblock Kluwer Academic Publishers, Dordrecht/Boston/London, 1988.

\bibitem{Hoch-Serre:Lie}
G.P.~Hochschild and J.-P. Serre.
\newblock{Cohomology of Lie Algebras}.
\newblock{\em Ann. of Math.}, 57:591--603, 1953.

\bibitem{Humphreys:Lie}
J.~E.~Humphreys.
\newblock{\em Introduction to Lie Algebras and Representation Theory}.
\newblock{\em Graduate Texts in Mathematics 9}.
\newblock{Springer-Verlag New York-Heidelberg-Berlin} 1972.

\bibitem{KS:Deformations}
K.~Kodaira and D.~C.~Spencer.
\newblock{On deformations of complex analytic structures, I-II; III}
\newblock{\em Ann. of Math.}, 67:328--466  1958; 71:43--76 1960.

\bibitem{Kostant:Borel}
B.~Kostant.
\newblock {Lie algebra cohomology and the generalized Borel-Weil theorem}.
\newblock {\em Ann. of Math.}, 74:329--387, 1961.

\bibitem{Mitchell:Categories}
B.~Mitchell.
\newblock{\em Theory of Categories}. 
\newblock{\em Pure and Applied Mathematics, Vol. XVII }
\newblock{Academic Press, New York-London} 1965.

\bibitem{Teich:quasikonforme}
O.~Teichm\"uller.
\newblock{Extremale quasikonforme Abbildungen und quadratische Differentiale}.
\newblock{\em Abh. Preuss. Akad. Wiss.}, 22:1--197, 1940.

\bibitem{Viviani:Melikian}
F.~Viviani
\newblock{Deformations of the restricted Melikian Lie algebra}.
\newblock{\em Comm. Algebra} 37:1850--1872, 2009.

\end{thebibliography}
\end{document}